\documentclass[journal,onecolumn,11pt]{IEEEtran}
\usepackage[utf8]{inputenc}
\usepackage[T1]{fontenc}
\usepackage{amsmath}
\usepackage{amssymb}
\usepackage{amsfonts}
\usepackage{graphicx}
\usepackage{hyperref}
\usepackage{geometry}

\begin{document}

\title{Dynamic Combinatorics of Quantum Modes in an Inverted Fock Space}

\author{
    Isayeva E.A.\\
    \textit{Institute of Physics of Ministry of Science and Education of Azerbaijan, H.Javid av.33,} \\
     Baku, Azerbaijan AZ 1143\\
    \textit{eaisaeva@mail.ru, e.isayeva@physics.science.az}
}

\maketitle

\begin{abstract}
An inversive approach to the combinatorics of Fock space in quantum field theory is proposed. The traditional description of Bose-Einstein and Fermi-Dirac systems is based on distributing indistinguishable quanta (particles) over energy states (modes). The present model introduces a methodological inversion: identical wave modes of a fixed frequency are declared as invariant elements (``balls''), while discrete occupancy numbers (Fock numbers) act as ``bins'' (spatial axes). It is shown that the time evolution of a quantum system governed by elastic interactions represents a trajectory of geometric hops of modes across occupancy axes. This approach eliminates ultraviolet divergences without introducing artificial regulators, visually differentiates thermodynamic quantum statistics, and rigorously derives Planck's entropy and the spectral characteristics of coherent laser systems.
\end{abstract}


\section{Introduction and the Conceptual Crisis of Quantization}
Modern approaches to the quantization of fields and spacetime, including Loop Quantum Gravity (LQG) and Causal Dynamical Triangulation (CDT), face fundamental conceptual barriers. Chief among them are the static ``problem of time'' (the disappearance of time evolution in the Wheeler-DeWitt equation) and uncontrollable ultraviolet divergences when attempting to discretize the continuum.

The root of these problems lies in the preservation of the classical paradigm, where space acts as an arena, and particles (or quanta of area/volume) are objects distributed across this arena. Historical analysis shows that Max Planck in 1900, introducing the concept of microscopic ``resonators'' (vibrators), modeled them as material elements of matter possessing discrete portions of energy $E = nh\nu$. The transition to the second quantization of fields (V. Fock, 1932) transferred this apparatus to the electromagnetic field itself, introducing the space of occupancy numbers. However, the mathematical exposition still operates with abstract creation and annihilation operators, masking the clear combinatorial structure of the underlying processes.

This paper develops an inversive model of quantum geometry in which quality (physical parameters of the mode) and quantity (the Fock number) swap places. The state space is generated dynamically in acts of elastic interaction, returning to quantum mechanics the transparency of the classical statistical mechanics of L. Boltzmann and J. W. Gibbs.

\section{Postulates of the Inversive Model}
The proposed theory is based on four fundamental postulates:
\begin{enumerate}
    \item \textbf{Invariance of Modes (The Mode as a Quantum Brick):} The elementary indistinguishable ``ball'' of the combinatorial system is a spacetime mode of a field with a fixed frequency $\nu$. The investigated physical system is an ensemble of $P$ such identical modes.
    \item \textbf{Topology of Occupancy Axes:} The role of geometric axes of the state space is played by quantization levels — fixed Fock numbers ($n = 0, 1, 2, \dots, N-1$). The coordinate of a particular mode at any given time is the discrete number of the axis (bin) where it is localized.
    \item \textbf{Immanent Discreteness of Time:} Time in the model is not an external continuous coordinate but is generated by a causal sequence of events. A single act of mode interaction is recognized as the unit of quantum time.
    \item \textbf{Elastic Interaction (Law of Invariant Conservation):} The elementary act of interaction represents a resonant (elastic) exchange of quanta between modes of the same frequency without the participation of external mediators. The total number of modes $P$ and the total number of quanta (the system's total energy $E_{\text{total}}$) remain strict invariants. The interaction shifts the two colliding modes to adjacent axes: one mode shifts to the $n-1$ axis, the second to the $n+1$ axis.
\end{enumerate}

\section{Combinatorial Apparatus and Separation of Quantum Statistics}
The mathematical apparatus of the model is based on an isomorphic inversion of the ``stars and bars'' method. For an ensemble of $P$ mode-balls distributed over $N$ bins (occupancy numbers), exactly $N-1$ identical configuration barriers (bars) are required.

The total number of elements in a linear sequence is $P + N - 1$. Since all modes are identical to each other and the bars are indistinguishable, the total number of unique macrostates (the combinatorial weight of the system $W$) is given by the formula of permutations with repetitions of two types of elements:
\begin{equation}
W = \frac{(P + N - 1)!}{P! \cdot (N - 1)!}
\end{equation}

The fundamental difference between quantum fields of matter and radiation is reduced to a simple combinatorial restriction on the geometric capacity of the occupancy axes:

\subsection{Fermi Structure (Fields of Matter / Electrons)}
Due to the Pauli exclusion principle, the occupancy bins have a rigid capacity limit. Mode-balls can occupy exactly two axes: $n=0$ (the vacuum state of the mode) or $n=1$ (the excited state). Repetition of elements on the same occupancy axis above unity is prohibited.

The combinatorics of the system narrows down to combinations without repetitions, expressed through permutations over two classes (occupied/free):
\begin{equation}
W_{\text{Fermi}} = \frac{P!}{m!(P-m)!}
\end{equation}
where $m$ is the number of occupied electronic modes. This strictly yields the Fermi-Dirac statistics.

\subsection{Bose Structure (Fields of Radiation / Photons)}
Planckian vibrators (modes of the electromagnetic field) possess an infinite capacity of occupancy axes. Any number of mode-balls from $0$ to $P$ can be placed in a single bin. The combinatorics of the system unfolds into full combinations with repetitions, described by expression (1), which corresponds to Bose-Einstein statistics.

\section{Thermodynamics of the Mode Ensemble and Derivation of Planck's Entropy}
The connection between the combinatorial weight of the inverted system and the macroscopic entropy $S$ is given by Boltzmann's postulate:
\begin{equation}
S = k_B \ln W
\end{equation}

Let us consider a thermodynamic ensemble of a Bose system, where the scale of numbers $P$ and $N$ allows us to neglect the unit in expression (1): $W \approx \frac{(P + N)!}{P! \cdot N!}$. Substituting the combinatorial weight into (3), we get:
\begin{equation}
S = k_B \big[ \ln(P + N)! - \ln(P!) - \ln(N!) \big]
\end{equation}

Applying Stirling's asymptotic approximation $\ln(X!) \approx X \ln X - X$ for large numbers, we expand the terms:
\begin{equation}
S \approx k_B \big[ (P + N) \ln(P + N) - (P + N) - (P \ln P - P) - (N \ln N - N) \big]
\end{equation}

The linear terms completely cancel out, reducing the expression to pure logarithmic geometry:
\begin{equation}
S = k_B \big[ (P + N) \ln(P + N) - P \ln P - N \ln N \big]
\end{equation}

Factoring the invariant number of modes $P$ out:
\begin{equation}
S = k_B \cdot P \left[ \left(1 + \frac{N}{P}\right) \ln(P + N) - \ln P - \frac{N}{P} \ln N \right]
\end{equation}

Using the identity $\ln(P+N) - \ln P = \ln(1 + \frac{N}{P})$, we finally transform the entropy equation of the quantum field to its canonical form:
\begin{equation}
S = k_B \cdot P \left[ \left(1 + \frac{N}{P}\right) \ln\left(1 + \frac{N}{P}\right) - \frac{N}{P} \ln\left(\frac{N}{P}\right) \right]
\end{equation}

The physical meaning of the ratio of the number of quantum bins to the number of modes $\frac{N}{P}$ in our inversion is strictly identical to the average occupancy number (Fock number) of a single mode at a given temperature: $\bar{n} = \frac{N}{P}$. Replacing the variable, we obtain the fundamental entropy formula derived by Max Planck in 1900:
\begin{equation}
S = k_B \cdot P \big[ (1 + \bar{n})\ln(1 + \bar{n}) - \bar{n}\ln\bar{n} \big]
\end{equation}

\section{Dynamic Equilibrium and the Radiation Law}
In the process of a cascade of continuous elastic hops of mode-balls across occupancy axes, the system evolves towards the most probable geometric distribution that maximizes entropy (9). By the definition of thermodynamic equilibrium, at a fixed temperature $T$, the relation holds:
\begin{equation}
\frac{1}{T} = \frac{dS}{dE_{\text{total}}}
\end{equation}

Given that the internal energy of the mode pool is $E = \bar{n} \cdot h\nu$, we transition to differentiation with respect to the average occupancy number ($dE = h\nu \cdot d\bar{n}$):
\begin{equation}
\frac{dS}{d\bar{n}} = k_B \cdot P \left[ \ln(1 + \bar{n}) + \frac{1+\bar{n}}{1+\bar{n}} - \ln\bar{n} - \frac{\bar{n}}{\bar{n}} \right] = k_B \cdot P \ln\left(\frac{1 + \bar{n}}{\bar{n}}\right)
\end{equation}

Substituting (11) into the equilibrium criterion (10), we get:
\begin{equation}
\frac{1}{T} = \frac{1}{h\nu} \frac{dS}{d\bar{n}} = \frac{k_B}{h\nu} \ln\left(\frac{1 + \bar{n}}{\bar{n}}\right)
\end{equation}

Solving equation (12) for the coordinate mark $\bar{n}$, we find the point of dynamic balance of the occupancy axes:
\begin{equation}
\ln\left(\frac{1}{\bar{n}} + 1\right) = \frac{h\nu}{k_B T} \implies \frac{1}{\bar{n}} + 1 = e^{\frac{h\nu}{k_B T}}
\end{equation}
\begin{equation}
\bar{n} = \frac{1}{e^{\frac{h\nu}{k_B T}} - 1}
\end{equation}

Formula (14) represents Planck's radiation law. The model proves that the macroscopic radiation spectrum is a direct consequence of the geometry of random walks of modes across the inverted axes of Fock space.
\section{Limiting Regimes and Resolution of the Ultraviolet Catastrophe}
The inversive combinatorial model automatically resolves the paradoxes of classical thermodynamics at extreme physical poles.
\subsection{Thermal Cryogenic Collapse ($T \rightarrow 0$)}
As the temperature approaches absolute zero, the total energy of the system is exhausted ($E_{\text{total}} \rightarrow 0$). All $P$ mode-balls are forcibly pulled into the initial vacuum bin $n=0$. The bars align in a monotonic row after the array of balls. The number of available arrangement variants becomes equal to one: $W \rightarrow 1$. Accordingly, the entropy of the system $S = k_B \ln(1) = 0$, which strictly satisfies the third law of thermodynamics (Nernst's theorem). Hops along the axes are completely blocked by the lack of energy resource.
\subsection{Resolution of the Ultraviolet Catastrophe ($\nu \rightarrow \infty$)}
In classical physics (the Rayleigh-Jeans law), it was assumed that high-frequency modes possess infinite heat capacity. In our model, as $\nu \rightarrow \infty$, the price of moving even a single mode-ball across the first bar from bin $n=0$ to bin $n=1$ requires an enormous quantum energy $h\nu$.
Since the total thermal energy of the system at a fixed temperature is limited, the barrier between the zero and first Fock numbers turns into a topologically insurmountable infinite wall. All high-frequency modes find themselves rigidly locked on the vacuum axis $n=0$. The system combinatorially ``switches off'' ultraviolet axes from the mixing process:
\begin{equation}
\lim_{\nu \rightarrow \infty} \bar{n} = \lim_{\nu \rightarrow \infty} \frac{1}{e^{\frac{h\nu}{k_B T}} - 1} = 0
\end{equation}
The radiation energy at high frequencies vanishes, eliminating the ultraviolet catastrophe at a purely geometric level.
\section{Application of the Model: Coherent Quantum Riot of Laser Systems}
Unlike the radiation of an absolute black body, where thermal mixing involves all frequencies through material mediators (cavity walls) and maximizes chaos ($W_{\text{BB}} = \prod W_\nu \rightarrow \infty$), laser systems implement a coherent ``quantum riot''.
In a laser cavity, external pumping breaks the statistical distribution and artificially isolates a single operating mode ($P=1$). The mechanism of stimulated emission avalanche-transports this single mode-ball through a cascade of bars deep to the right — to an astronomically high occupancy axis $N$ (on the order of $10^9$ quanta or more).
Let us substitute the laser radiation parameters ($P=1$, number of bins $N$) into the basic combinatorial equation (1):
\begin{equation}
W_{\text{laser}} = \frac{(1 + N - 1)!}{1! \cdot (N - 1)!} = \frac{N!}{(N - 1)!} = N
\end{equation}
The number of possible microstates of a laser beam is strictly equal to the number of generated photons $N$. Let us calculate the thermodynamic entropy of the resulting macrostate:
\begin{equation}
S_{\text{laser}} = k_B \ln N
\end{equation}
For a macroscopic laser beam with a density of $N = 10^9$ photons, the value of the logarithm is $\ln(10^9) \approx 20.7$. Multiplying by Boltzmann's constant $k_B \approx 1.38 \cdot 10^{-23}$ J/K, we get:
\begin{equation}
S_{\text{laser}} \approx 2.85 \cdot 10^{-22} \text{ J/K} \approx 0
\end{equation}
The entropy of laser radiation drops to physical zero. A laser beam is not a chaotic conglomerate of particles, but a single spatial mode, monopolistically fixed at a giant occupancy coordinate. This combinatorially explains the ideal coherence and absence of thermodynamic noise in laser radiation (Glauber states).
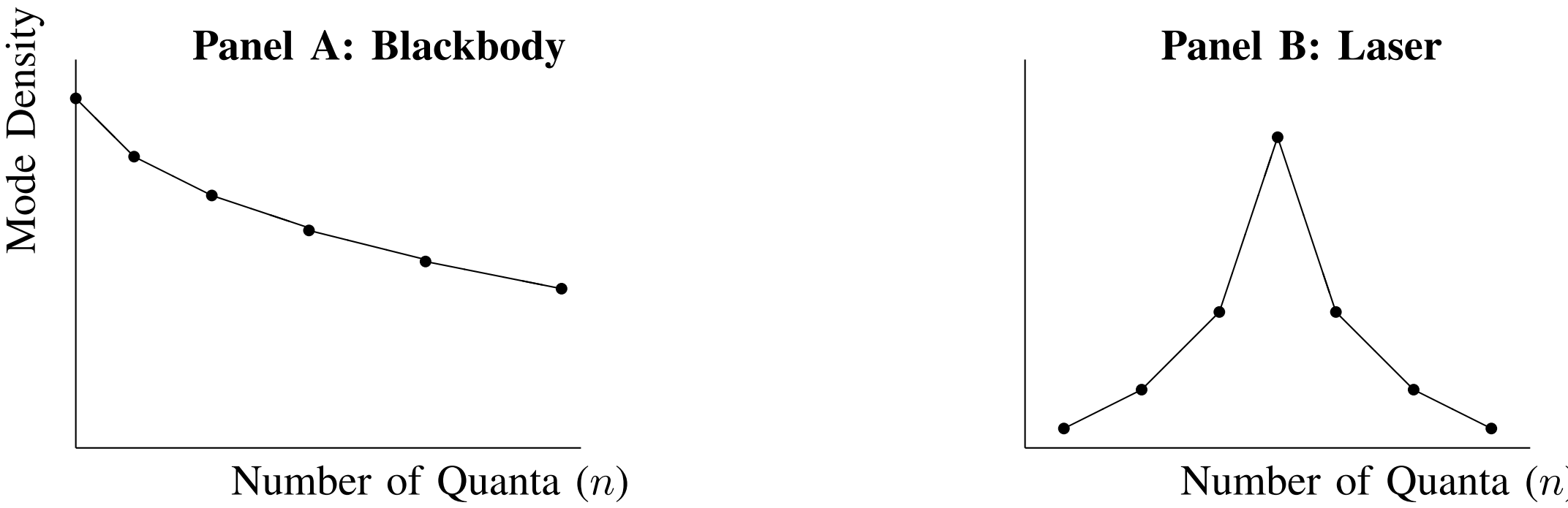
\begin{figure}[htbp]
\centering
\begin{minipage}{0.48\textwidth}
    \centering
    % Панель А: Экспоненциальное тепловое затухание
    \begin{picture}(150,120)
        \put(10,10){\line(1,0){130}} % ось X
        \put(10,10){\line(0,1){100}} % ось Y
        \put(10,100){\circle*{3}} \put(10,100){\line(1,-1){15}}
        \put(25,85){\circle*{3}}  \put(25,85){\line(2,-1){20}}
        \put(45,75){\circle*{3}}  \put(45,75){\line(3,-1){25}}
        \put(70,66){\circle*{3}}  \put(70,66){\line(4,-1){30}}
        \put(100,58){\circle*{3}} \put(100,58){\line(5,-1){35}}
        \put(135,51){\circle*{3}}
        \put(50,-2){Number of Quanta ($n$)}
        \put(-8,60){\rotatebox{90}{Mode Density}}
        \put(40,110){\textbf{Panel A: Blackbody}}
    \end{picture}
\end{minipage}
\hfill
\begin{minipage}{0.48\textwidth}
    \centering
    % Панель Б: Когерентный пик лазера (Пуассон)
    \begin{picture}(150,120)
        \put(10,10){\line(1,0){130}} % ось X
        \put(10,10){\line(0,1){100}} % ось Y
        \put(20,15){\circle*{3}}  \put(20,15){\line(2,1){20}}
        \put(40,25){\circle*{3}}  \put(40,25){\line(1,1){20}}
        \put(60,45){\circle*{3}}  \put(60,45){\line(1,3){15}}
        \put(75,90){\circle*{3}}  \put(75,90){\line(1,-3){15}}
        \put(90,45){\circle*{3}}  \put(90,45){\line(1,-1){20}}
        \put(110,25){\circle*{3}} \put(110,25){\line(2,-1){20}}
        \put(130,15){\circle*{3}}
        \put(50,-2){Number of Quanta ($n$)}
        \put(45,110){\textbf{Panel B: Laser}}
    \end{picture}
\end{minipage}
\caption{Distribution of the mode ensemble over occupancy axes (Fock numbers) in the inverted Fock space: Panel A shows the thermal exponential decay characteristic of blackbody radiation (maximum entropy state), while Panel B demonstrates the concentrated Poisson-like distribution of a single isolated mode in a coherent laser system (minimum entropy state).}
\label{fig:modes_distribution}
\end{figure}
\section{Conclusions and Conceptual Advantages over Quantum Gravity Models}
The proposed inversive model of mode combinatorics has three key advantages over leading modern quantum theories:
\begin{enumerate}
\item \textbf{Resolution of the ``problem of time'':} Unlike Loop Quantum Gravity, which completely excludes time from fundamental equations, the present model retains the arrow of time as an ordered sequence of elastic interactions of modes. Time here is immanently linked to the geometric motion of the system.
\item \textbf{Occam's Razor against mathematical overload:} The theory replaces the cumbersome apparatus of non-commutative geometry of creation/annihilation operators with a transparent isomorphism of permutations with repetitions. The properties of quantum ensembles are derived from the basic principles of probability theory of Boltzmann and Gibbs.
\item \textbf{Natural absence of singularities:} The separation of quality (modes) and quantity (Fock numbers) blocks the appearance of mathematical infinities. Mode interactions themselves balance the capacity of the state space, forming smooth macroscopic distributions without invoking artificial regularization methods.
\end{enumerate}
Inversive quantum combinatorics of modes opens up a promising vector of research in the field of dynamic causal networks, proving that the topological dimensionality and physical properties of our Universe are capable of weaving themselves from simple acts of shuffling wave modes along discrete occupancy axes.

\end{document}